\input amstex

\documentstyle{amsppt}
\document
\topmatter
\title
Compact conformally K\"ahler Einstein-Weyl manifolds
\endtitle
\author
W\l odzimierz Jelonek
\endauthor
\abstract{We give a description of compact conformally K\"ahler
Einstein-Weyl manifolds whose Ricci tensor is Hermitian.}
 \endabstract
\thanks{The paper was supported by Narodowe Centrum Nauki grant
DEC-2011/01/B/ST1/02643}\endthanks
\endtopmatter
\define\G{\Gamma}
\define\e{\epsilon}

\define\n{\nabla}

\define\w{\wedge}
\define\k{\diamondsuit}
\define\th{\theta}
\define\p{\parallel}
\define\a{\alpha}

\define\m{(M,g,J)}
\define\om{\omega}
\define\lb{\lambda}

\define\De{\Cal D}
\define\DE{\Cal D^{\perp}}

\define\AC{\Cal A\oplus\Cal C^{\perp}}
\define\1{D_{\lb}}
\define\2{D_{\mu}}
\define\Om{\Omega}

\define\0{\omega}

\define\Lb{\Lambda}

\medskip
{\bf 0. Introduction.} In this paper we shall investigate compact
Einstein-Weyl structures $(M,[g],D)$ on a complex manifold
$(M,J)$, $dim M\ge4$,  which are conformally K\"ahler and whose
Ricci tensor $\rho^D$ is Hermitian i.e. $\rho^D$ is $J$
-invariant.

We give a complete classification of  compact Einstein-Weyl
structures $(M,[g],D)$ with dim $M\ge4$ such that $(M,[g],J)$ is
conformally K\"ahler, i.e. there exists a metric $g_0\in [g]$ such
that $(M,g_0,J)$ is K\"ahler and whose Ricci tensor $\rho^D$ is
$J$-invariant  i.e. $$\rho^D(JX,JY)=\rho^D(X,Y).$$  Conformally
K\"ahler Einstein manifolds were classified by A. Derdzi\'nski and
G.  Maschler in [D-M-3]. Compact Gray bi-Hermitian manifolds are
partially classified in [J-2]. It is proved in [J-1]  that compact
Einstein-Weyl manifolds are also Gray manifolds (see [Gr]). The
compact Einstein-Weyl 3-dimensional manifolds are studied in [T].
The compact Einstein-Weyl manifolds on complex manifolds
compatible with complex structure are studied in  [PS1], [PS2],
[MPPS], [WW]. In [D-T] there are studied Riemannian manifolds $(M,
g)$ in four dimensions which are locally conformally  K\"ahler.
The Einstein-Weyl structures which are conformally K\"ahler
studied in [MPPS], [WW]
 have the Ricci tensor $\rho^D$ of the Weyl structure $(M,[g],D)$
which $J$ -invariant. In the first section of the paper we recall
some facts from [J-1] and describe compact Einstein-Weyl manifolds
with the Gauduchon metric as Gray manifolds. In section 2 we
describe the Riemannian structure of Gray manifold corresponding
to Einstein-Weyl manifold with the Gauduchon metric. In section 3
we prove that compact Einstein-Weyl manifold with Hermitian Ricci
tensor admits a holomorphic Killing vector field with special
K\"ahler-Ricci potential and consequently $M=\Bbb{CP}^n$ or
$M=\Bbb{P}(L\oplus\Cal O)$ where $L$ is a holomorphic line bundle
over K\"ahler-Einstein manifold.  In section 4 we prove that if
$M=\Bbb{P}(L\oplus\Cal O)$ then $L$ is a holomorphic line bundle
over K\"ahler-Einstein manifold with positive scalar curvature.

 \par
 \medskip
  {\bf 1.  Einstein-Weyl geometry and Killing tensors.} We start with some basic facts
  concerning Einstein-Weyl geometry. For more details see [T], [PS1],
[PS2].

Let $M$ be a $n$-dimensional manifold with a conformal structure
$[g]$ and a torsion-free affine
 connection $D$. This defines an Einstein-Weyl (E-W) structure
if  $D$ preserves the conformal structure i.e. there exists a 1-form $\0$ on $M$
such that
$$Dg=\0\otimes g \tag 1.1$$
and the Ricci tensor $\rho^D$ of $D$ satisfies the condition
$$\rho^D(X,Y)+\rho^D(Y,X)=\bar\Lb g(X,Y) \text{  for every  } X,Y\in TM$$
for some function $\bar\Lb\in C^{\infty}(M)$. P.Gauduchon proved
([G]) the fundamental theorem that if $M$ is compact then there
exists a Riemannian metric $g_0\in [g]$ for which $\delta \0_0=0$
and $g_0$ is unique up to homothety.
 We shall call $g_0$ a standard metric of E-W structure $(M,[g],D)$. Let $\rho$ be the Ricci tensor of $(M,g)$
  and let us denote by $S$ the Ricci
 endomorphism of $(M,g)$, i.e. $\rho(X,Y)=g(X,SY)$.
We recall two important theorems (see [T], [PS1]):
\par
\medskip
{\bf Theorem 1.1. }  {\it A metric g and a 1-form $\0$ determine
an E-W structure if and only if there exists a function $\Lb\in
C^{\infty}(M)$ such that
$$\rho^{\n}+\frac14(n-2) \Cal D\0=\Lb g  \tag 1.2$$
where $\Cal D\0(X,Y)=(\n_X\0)Y+(\n_Y\0)X+\0(X)\0(Y)$ and $n=dim \
M$. If (1.2) holds then}
$$\bar\Lb=2\Lb+div\0-\frac12 (n-2)\p\0^{\sharp}\p^2\tag 1.3$$
\par
\medskip
 By
$\tau=\operatorname{tr}_g\rho$ we shall denote the scalar
curvature of $(M,g)$.  Compact Einstein-Weyl manifolds with the
Gauduchon metric are Gray manifolds.
\medskip
{\it Definition.}  A Riemannian manifold $(M,g)$ will be  called a
Gray $\Cal {A}\oplus\Cal{C}^{\perp}$ manifold if the tensor
$\rho-\frac{2\tau}{n+2}g$ is a Killing tensor.
\medskip
{\it Definition.} A self-adjoint $(1,1)$ tensor on a Riemannian
manifold $(M,g)$ is called a Killing tensor if $$g(\n
S(X,X),X)=0$$ for arbitrary $X\in TM$.

The condition  $g(\n S(X,X),X)=0$ is equivalent to $$\frak
C_{X,Y,Z}g(\n S(X,Y),Z)=0$$ for arbitrary $X,Y,Z\in \frak X(M)$,
where $\frak C$ denotes the cyclic sum.

In this paper Gray $\Cal {A}\oplus\Cal{C}^{\perp}$ manifolds will
be  called for short Gray manifolds or $\Cal
{A}\oplus\Cal{C}^{\perp}$ manifolds.  Gray manifolds were first
defined by Gray ([Gr]).
\medskip
K.P. Tod proved [T] that the Gauduchon metric admits a Killing
vector field, more precisely he proved:.
\par
\medskip
{\bf Theorem 1.2. } {\it Let $M$ be a compact E-W manifold and let
$g$ be the standard metric with the corresponding 1-form $\0$.
Then the vector  field $\0^{\sharp}$ dual to the form $\om$ is
 a Killing vector field on $M$.}
\par
\medskip
   In the next two theorems, which are proved also in [J-1],  we characterize
    compact Einstein-Weyl
   manifolds $(M,g)$ with the Gauduchon metric as Gray ([Gr]) manifolds
   and show that the eigenvalues $\lb_0,\lb_1$ of the Ricci tensor satisfy the equation
   $(n-4)\lb_1+2\lb_0=C_0=const$ which we shall use later. We sketch the proofs of the theorems.
    Our  motivation is to use the structure of Gray manifolds to use
   ideas from [J-2], where we classified a different class of Gray
   manifolds.
        From the above theorems it follows  (see [J-1])
\medskip
{\bf Theorem 1.3. } {\it Let $(M,[g])$ be a compact E-W manifold,
$n=$dim$M\ge3$, and let $g$ be the standard metric on $M$. Then
$(M,g)$ is an $\AC$-manifold. The manifold $(M,g)$ is Einstein or
the Ricci tensor $\rho^{\n}$ of $(M,g)$ has exactly two
eigenfunctions
 $\lb_0\in C^{\infty}(M),\lb_1=\Lb$ satisfying the following conditions:

(a) $(n-4)\lb_1+2\lb_0=C_0=const$

(b) $\lb_0\o\lb_1$  on $M$

(c) $dim\ ker(S-\lb_0Id)=1,dim\ ker(S-\lb_1Id)=n-1$ on
$U=\{x:\lb_0(x)\ne\lb_1(x)\}$,

(d) $\lb_1-\lb_0=\frac{n-2}4\p\xi\p^2$ where
$\xi=\om^{\sharp}\in\frak{iso}(M)$.
\medskip
In the addition $\lb_0=\frac1n Scal^D_g$ where $Scal^D_g=tr_g\rho^D$
 denotes the conformal scalar curvature of $(M,g,D)$.}
\par
\medskip
 {\it Proof.}  Note that $\0(X)=g(\xi,X)$ where $\xi\in\frak{iso}(M)$ and the formula
$$\rho^{\n}+\frac14(n-2)\0\otimes\0=\Lb g\tag 1.4$$
holds.   Thus $\n_X(\0\otimes\0)(X,X)= 0$. From (1.4) it follows
that
$$\n_X\rho(X,X)=X\Lb g(X,X). \tag 1.5$$
It means that $(M,g)\in\AC$ and $d(\Lb-\frac2{n+2}\tau)=0$ ([J-1],
Lemma 1.5) , where $\tau$ is the scalar curvature of $(M,g)$.
From (1.5) it follows that the tensor $T=S-\Lb Id$ is a Killing
tensor. Let us denote by $\xi$ the Killing vector field dual to
$\0$. Note that
$\rho(\xi,\xi)=(\Lb-\frac14(n-2)\p\xi\p^2)\p\xi\p^2$ and if
$X\perp\xi$ then $SX=\Lb X$. Hence the tensor $S$ has two
eigenfunctions $\lb_0=\Lb-\frac14(n-2)\p\xi\p^2$ and $\lb_1=\Lb$.
This proves (b).  Note that
$$\tau=\lb_0+(n-1)\lb_1=n\Lb-\frac14(n-2)\p\xi\p^2.$$
and $2\tau-(n+2)\Lb=C_0=const$. Thus
$C_0=(n-2)\Lb-\frac12(n-2)\p\xi\p^2$. However
$(n-4)\lb_1+2\lb_0=(n-2)\Lb-\frac12(n-2)\p\xi\p^2$ which proves
(a). Note also that
$$\frac1ns^D_g=\Lb-\frac{n-2}4\p\xi\p^2=\lb_0\tag 1.6$$
which finishes the proof.$\k$

On the other hand the following theorem holds (see [J-1]).
\par
\medskip
{\bf Theorem 1.4. } {\it  Let $(M,g)$ be a compact $\AC$ manifold.
Let us assume that the Ricci tensor $\rho$ of $(M,g)$ has exactly
two eigenfunctions
 $\lb_0,\lb_1$ satisfying the conditions:

(a) $(n-4)\lb_1+2\lb_0=C_0=const$

(b) $\lb_0\o\lb_1$  on $M$

(c) $dim\ ker(S-\lb_0Id)=1,dim\ ker(S-\lb_1Id)=n-1$ on
$U=\{x:\lb_0(x)\ne\lb_1(x)\}$.

Then there exists a two-fold Riemannian covering $(M',g')$ of
$(M,g)$ and a Killing vector field $\xi\in\frak{iso}(M')$ such
that $(M',[g'])$ admits two different E-W structures with the
standard metric $g'$ and the corresponding 1-forms
$\0_{\mp}=\mp\xi^{\flat}$ dual to the vector fields $\mp\xi$.
Additionally $\lb_1-\lb_0=\frac{n-2}4\p\xi\p^2$. The condition (b)
may be replaced by the condition

(b1) there exists a point $x_0\in M$ such that $\lb_0(x_0)<\lb_1(x_0)$.}
\par
\medskip
{\it Proof.} Note that $\tau = (n-1)\lb_1+\lb_0$ and
$C_0=(n-4)\lb_1+2\lb_0$. It follows that
$$\lb_1=\frac{2\tau-C_0}{n+2},\ \ \ \lb_0=\frac{(n-1)C_0-(n-4)\tau}{n+2}.
\tag 1.7$$ In particular $\lb_0,\lb_1\in C^{\infty}(M)$. Let $S$
be the Ricci endomorphism of $(M,g)$ and let us define the tensor
$T:=S-\lb_1Id$. Since from (1.7) we have $d\lb_1=\frac2{n+2}d\tau$
it follows that $T$ is a Killing tensor with two eigenfunctions:
$\mu = 0$ and $\lb = \lb_0 - \lb_1$. Hence  there exists a two
fold Riemannian covering $p:(M',g')\rightarrow (M,g)$ and a
Killing vector field $\xi\in\frak{iso}(M')$  (see [J-1], th. 2.10)
such that $S'\xi=(\lb_0\circ p)\xi$ where $S'$ is the Ricci
endomorphism of $(M',g')$. Note also that $\p\xi\p^2 = |\lb-\mu| =
|\lb_0-\lb_1|$. Let us define the 1-form $\0$ on $M'$ by
$\0=c\xi^{\flat}$ where $c=2\sqrt{\frac1{n-2}}$. It is easy to
check that with such a choice of $\0$ equation (1.4) is satisfied
and $\delta\0=0$. Thus $(M',g',\0)$ defines an E-W structure and
$g'$ is the standard metric for $(M',[g'])$. Note that
$(M,g',-\0)$ gives another E-W structure corresponding to the
field $-\xi$. $\k$
\par
\medskip
{\it Corollary 1.5. }  Let $(M,g)$ be a compact simply connected
manifold satisfying the assumptions of Th.1.4. Then $(M,[g])$
admits two E-W structures with the standard metric $g$.

\par
\medskip
{\bf 2. Killing tensors.}  In this section we describe the
Riemannian manifold $(M,g)$ where $g\in[g])$ is the standard
metric of Einstein-Weyl structure $(M,[g],D)$.  We describe such
manifolds under an additional condition, that the distribution
corresponding to the eigenvalue $\lb_1$ of the Ricci tensor is
totally geodesic.   Note that if $(M,g)$ is a compact
Einstein-Weyl manifold with the Gauduchon metric then the
distribution $\De_{\lb_1}$ is totally geodesic since is orthogonal
to the distribution spanned by a Killing vector field.

We say, that a distribution (not necessarily  integrable)  $\Cal
D$ is totally geodesic, if $\n_XX\in\G(\Cal D)$ for every
$X\in\G(\Cal D)$.

We start with:
\par
\medskip
{\bf Lemma 2.1.} {\it Let S be a self-adjoint tensor on }$(M,g)$
{\it with exactly two eigenvalues} $\lb,\mu$. {\it If the
distributions } $\Cal D_{\lb},\Cal D_{\mu}$ {\it are both
umbilical, $\n\lb\in\G(\2),\n \mu\in\G(\1)$ and the mean curvature
normals  $\xi_{\lb},\xi_{\mu}$ of the distributions $\1,\2$
respectively satisfy the equations}
$$\xi_{\lb}=\frac1{2(\mu-\lb)}\n\lb,\qquad\xi_{\mu}=\frac1{2(\lb-\mu)}\n\mu,$$
{\it then $S$ is a Killing tensor.}
\medskip
{\it Proof.}  We have to show that $g(\n S(Z,Z),Z)=0$ for
arbitrary $Z\in TM$. Let $Z=X+Y$ where $X\in\Cal D_{\lb},Y\in\Cal
D_{\mu}$. Then
$$\gather g(\n S(Z,Z),Z)=g(\n S(X,X),X)+2g(\n S(X,X),Y)+g(\n S(Y,X),X)\\+
2g(\n S(Y,Y),X)+g(\n S(X,Y),Y)+g(\n S(Y,Y),Y).\endgather$$ Since
$\n S(X,X)=(\lb-\mu)g(X,X)\xi_{\lb},\n
S(Y,Y)=(\mu-\lb)g(Y,Y)\xi_{\mu}$ and $$g(\n S(X,X),X)=0, g(\n
S(Y,Y),Y)=0$$ one can easily check that $g(\n S(Z,Z),Z)=0$.
\medskip
{\it Remark.} We call here a vector field $\xi$ the mean curvature
normal of a umibilical distribution $\Cal D$ if for every
$X\in\G(\Cal D)$ we have $\pi(\n_XX)=g(X,X)\xi$ where $\pi$ is a
projection on the orthogonal complement of $\Cal D$. Note that
$\Cal D$ may not be integrable.

\medskip

{\bf Proposition 2.2. } {\it Let $(M,g)$ be a $2n$-dimensional
Riemannian manifold whose Ricci tensor $\rho$ has two eigenvalues
$\lb_0(x),\lb_1(x)$ of   multiplicity 1 and $2n-1$ respectively at
every point $x$ of $M$. Assume that the eigendistribution
 $\De_{\lb_1}$ corresponding to $\lb_1$ is
totally geodesic. Then  $(M,g)$ is a Gray manifold if and only if
$2\lb_0+(2n-4)\lb_1$ is constant and $\n\tau\in\G(\De_{\lb_1})$. }
\bigskip
{\it Proof. }  Let $S_0$ be the Ricci endomorphism of $(M,g)$,
i.e. $\rho(X,Y)=g(S_0X,Y)$. Let  $S$ be the tensor defined by the
formula
$$S_0=S+\frac{\tau}{n+1} \operatorname{id}.\tag 2.1$$
Then $$\operatorname{tr} S=-\frac{(n-1)\tau}{n+1}.\tag 2.2$$ Let
$\lb_0,\lb_1$ be the eigenfunctions of $S_0$ and let us assume
that
$$2\lb_0+(2 n-4)\lb_1=C\tag 2.3$$
where $C\in\Bbb{R}$. Note that $S$ also has two eigenfunctions
which we denote by $\lb_0',\lb_1'$ respectively. It is easy to see
that $\lb_0'=-\frac
{n-1}{n+1}\tau+C\frac{2n-1}{2(n+1)},\lb_1'=-\frac{C}{2(n+1)}$ and
$\lb_0=-\frac{\tau(n-2)}{n+1}+C\frac{2n-1}{2(n+1)},\lb_1=\frac{\tau}{(n+1)}-\frac{C}{2(n+1)}$.
Since the distribution $\De_{\lb_0}$ is umbilical we have
$\n_XX_{|\De_{\lb_1}}=g(X,X)\xi$ for any $X\in\G(\De_{\lb_0})$
where $\xi$ is the mean curvature normal of $\De_{\lb_0}$. Since
the distribution $\De_{\lb_1}$ is totally geodesic we also have
$\n_XX_{|\De_{\lb_0}}=0$ for any $X\in\G(\De_{\lb_1})$. Let
$\{E_1,E_2,E_3,E_4,...,E_{2n-1},E_{2n}\}$ be a local  orthonormal
basis of $TM$ such that $\De_{\lb_0}=\text{span }\{E_1\}$ and
$\De_{\lb_1}=\text{span }\{E_2,E_3,E_4,...,E_{2n}\}$. Then
$\n_{E_i}E_{i|\De_{\lb_0}}=0$ for $i\in\{2,3,4,...,2n\}$ and
$$\n_{E_1}E_{1|\De_{\lb_1}}=\xi.$$  Consequently (note that
$\n\lb'_{0|\De_{\lb_0}}=0$ if and only if
$\n\tau_{|\De_{\lb_0}}=0$),
$$\gather \operatorname{tr}_g\n S=\sum_{i=1}^n\n S(E_i,E_i)=-
(S-\lb_0
\operatorname{id})(\n_{E_1}E_{1})+\n\lb_{0|\De_{\lb_0}}\tag
2.4\\=-(\lb_1-\lb_0)\xi\endgather$$ if we assume that
$\n\tau_{|\De_{\lb_0}}=0$. On the other hand,
$\operatorname{tr}_g\n S_0=\frac{\n\tau}2$ and
$\operatorname{tr}_g\n S=\operatorname{tr}_g\n
S_0-\frac{\n\tau}{n+1}$. Consequently,
$$\operatorname{tr}_g\n S=\frac{(n-1)\n\tau}{2(n+1)}=-\frac12\n\lb_0'.\tag 2.5$$
Thus $\xi=-\frac1{2(\lb_0-\lb_1)}\n\lb_0'$. From  the Lemma it
follows that $(M,g)$ is an $\AC$-manifold if $2\lb_0+(n-4)\lb_1$
is constant and $\n\tau\in\G(\De_{\lb_1})$. These conditions are
also necessary since $\n\lb_1'=0$ if $(M,g)$ is an $\AC$-manifold
and $\De_{\lb_1}$ is totally geodesic. Analogously
$\xi=-\frac1{2(\lb_0'-\lb_1')}\n\lb_0'$ and $\n\lb_0'=-\frac {n-1}
{(n+1)}\n\tau\in\G(\De_{\lb_1})$, where $\xi$ is the mean
curvature normal of the umbilical distribution $\De_{\lb_0}$, if
$(M,g)$ is an $\AC$-manifold.$\k$

\par
\medskip
{\bf 3. Conformally K\"ahler Einstein-Weyl manifolds.}  Let $g$ be
the standard metric of $(M,[g])$. Now let us recall that
$\rho^D(X,Y)=\lb_0g(X,Y)+\frac n4d\0(X,Y)$. Let us assume that
$(M,J)$ is complex and $[g]$ is Hermitian i.e. $ g(JX,JY)=g(X,Y)$.
It follows that $\rho^D$ is $J$-invariant if and only if $d\0$ is
a $(1,1)$ form, $d\0(JX,JY)=d\0(X,Y)$.  Since $\0(X)=g(\xi,X)$ it
follows that $d\0$ is a $(1,1)$ form iff $\n_{JX}\xi=J\n_X\xi$.
\medskip
\par
{\bf Proposition 3.1.} {\it Let $(M,J)$ be a compact complex
manifold with conformal Hermitian structure $[g]$.  Let us assume
that $[g]$ is conformally K\"ahler and   $f^2g$ is a K\"ahler
metric on $(M,J)$ where $g$ is the standard metric and $f\in
C^{\infty}(M)$. If $(M,[g])$ is Einstein-Weyl with $J$-invariant
Ricci tensor $\rho^D$ then $J\xi$ is colinear with $\n f$ in
$U=\{x:\xi_x\ne 0\}$   and $\xi$ is a holomorphic Killing field on
$(M,f^2g,J)$.}

\medskip
{\it Proof.}  Let $\n$ be a Levi-Civita  connection of the
standard metric $g$ and $\n^1$ be a Levi-Civita connection of the
K\"ahler metric  $g_1=f^2g$.  Note that $\xi$ is a conformal field
on $(M,g_1)$, $L_{\xi}g_1=L_{\xi}(f^2g)=2\xi \ln  f g_1=\sigma
g_1$.  Every conformal field on a compact K\"ahler manifold is
Killing (see [L]), hence consequently $\xi f=0$ and
$\xi\in\frak{iso}(M,g_1)$.  On a K\"ahler compact manifold every
Killing vector field is holomorphic (see[Mor  ]).  Thus
$\xi\in\frak{hol}(M,J)$.  Note that

$$\n_X\xi=\n^1_X\xi -d\ln f(X)\xi -d\ln f(\xi)X+g_1(X,\xi)\n^1\ln
f.$$ Thus
$$\gather  \n_{JX}\xi-J\n_X\xi=-d\ln f(JX)\xi -d\ln f(\xi)JX+g_1(JX,\xi)\n^1\ln
f\\+d\ln f(X)J\xi +d\ln f(\xi)JX-g_1(X,\xi)J\n^1\ln f.\endgather$$
Hence  $\n_{JX}\xi=J\n_X\xi$ if $$-d\ln
f(JX)\xi+g_1(JX,\xi)\n^1\ln f+d\ln f(X)J\xi -g_1(X,\xi)J\n^1\ln
f=0.$$ Put $X=\xi$ then we get $g_1(\xi,\xi)J\n^1\ln f=-d\ln
f(J\xi)\xi$.  It follows that in $U=\{x\in M: \xi_x\ne 0\}$ there
exists a smooth function $\phi$ such that $\n^1f=\phi J\xi$.$\k$

Let us recall the definition of a special K\"ahler-Ricci potential
([D-M-1], [D-M-2]).
 \medskip {\it Definition.}  A nonconstant
function $\tau\in C^{\infty}(M)$, where $\m$ is a K\"ahler
manifold, is called a special K\"ahler-Ricci potential if the
field $X= J(\n \tau)$ is a Killing vector field and at every point
with $d\tau\ne 0$ all nonzero tangent vectors orthogonal to the
fields $X,JX$ are  eigenvectors of both $\n d\tau$ and the Ricci
tensor $\rho$ of $\m$.

\medskip
Now our aim is to prove
\medskip
{\bf Theorem 3.2} {\it Let us assume that $(M,[g],J)$ is a
compact, conformally K\"ahler Einstein-Weyl manifold with
Hermitian Ricci tensor $\rho^D$ which is not conformally Einstein.
Then the conformally equivalent K\"ahler manifold  $(M,g_1,J)$
admits a holomorphic Killing field with a K\"ahler-Ricci
potential.  Thus $M=\Bbb P(L\oplus\Cal O)$ where $L$ is a
holomorphic line bundle over a compact  K\"ahler Einstein manifold
$(N,h)$ of positive scalar curvature  or is a complex projective
space $\Bbb{CP}^n$.}
\medskip
{\it Proof.} Let $\rho,\rho^1$ be the Ricci tensors of conformally
related riemannian metrics $g,g_1=f^2g$.  Then

$$ \rho=\rho^1+(n-2)f^{-1}\n^1df+[f^{-1}\Delta^1f-(n-1)f^{-2}g_1(\n^1
f,\n^1 f)]g_1.$$

Note that for arbitrary $X,Y\in\frak X(M)$ we have
$\n^1df(X,Y)=g_1(\n^1_X\n^1f,Y)=g_1(X\phi J\xi,Y)+\phi
g_1(J\n^1_X\xi,Y)$.  Thus for any $X,Y\in\frak X(M)$:

$$\gather \rho(X,Y)-(n-2)fX\phi g(J\xi,Y)=\rho^1(X,Y)+(n-2)f^{-1}\phi
g_1(J\n^1_X\xi,Y)+\tag 3.1\\[f^{-1}\Delta^1f-(n-1)f^{-2}g_1(\n^1 f,\n^1
f)]g_1(X,Y),\endgather$$ where $\Delta^1f=tr_{g_1}\n^1df$.

We shall show that $\xi$ has zeros on  $M$.   If $\xi\ne0$ on $M$
then the function $\phi$ would be defined and smooth on the whole
of $M$.   Since $M$ is compact it would imply that there exists a
point $x_0\in M$ such that $d\phi=0$ at $x_0$.  On the other hand
the eigenvalues $\lb_0,\lb_1$ of the Ricci tensor $\rho$ satisfy
$\lb_0-\lb_1= Cg(\xi,\xi)$ where $C\ne 0$ is a real number.  Since
$\xi\ne0$ it follows that the eigenvalues of $\rho$ do not
coincide at any point of $M$.  In particular $\rho$ is not
$J$-invariant at $x_0$, a contradiction, since the right hand part
of (3.1) is $J$-invariant.  It implies that $\xi$ is a holomorphic
Killing vector field with zeros and thus has a potential $\tau$
(see [K]), i.e.  there exists $\tau\in C^{\infty}(M)$ such that
$\xi=J\n^1\tau$.  Hence $df=-\phi d\tau$ and $d\phi\w d\tau=0$. It
implies that $d\phi=\a d\tau$.  Thus we have for arbitrary
$X,Y\in\frak X(M)$:

$$\gather \rho(X,Y)+(n-2)f^{-1}\a d\tau(X)d\tau(Y)=\rho^1(X,Y)\tag 3.2\\-(n-2)f^{-1}\phi
H^{\tau}(X,Y)-[f^{-1}\a
Q+f^{-1}\phi\Delta^1\tau+(n-1)f^{-2}\phi^2Q]g_1(X,Y).\endgather$$
where $Q=g_1(\xi,\xi)$.

Note that the tensor
$\widetilde{\rho}(X,Y)=\rho(X,Y)+(n-2)f^{-1}\a d\tau(X)d\tau(Y)$
is $J$-invariant.  In particular
$\widetilde{\rho}(\xi,\xi)=\lb_0g(\xi,\xi)=\lb_o\frac{Q}{f^2}$. On
the other hand
$\widetilde{\rho}(\n^1\tau,\n^1\tau)=\lb_1\frac{Q}{f^2}+(n-2)f^{-1}\a
Q^2$.  Hence  $(\lb_0-\lb_1)\frac{Q}{f^2}=(n-2)f^{-1}\a Q^2$.
Since $\lb_0-\lb_1=-\frac14(n-2)\frac{Q}{f^2}$ we get
$\a=-\frac1{4f^3}$.  Hence
$$d\phi=-\frac1{4f^3}d\tau=\frac1{4f^3}\frac{df}{\phi},\tag 3.3$$
and we get  $8\phi d\phi=-d(\frac1{f^2})$.  Hence
$d(4\phi^2+\frac1{f^2})=0$ and $4 \phi^2+\frac1{f^2}=C=const$.

Let us denote $\chi=(n-2)f^{-1}\phi, \sigma_0=f^{-1}\a
Q+f^{-1}\phi\Delta^1\tau+(n-1)f^{-2}\phi^2Q$. Note also that the
vector field $v=\n^1\tau$ is holomorphic and consequently
$i_v\rho^1=-\frac12d\Delta^1\tau=-\frac12dZ$ where
$Z=\Delta^1\tau$.  From the equation
$$\gather \widetilde{\rho}(X,Y) =\rho^1(X,Y)-(n-2)f^{-1}\phi
H^{\tau}(X,Y)-\tag 3.4\\[f^{-1}\a Q+f^{-1}\phi\Delta^1\tau+(n-1)f^{-2}\phi^2Q]g_1(X,Y).\endgather$$

valid for arbitrary $X,Y\in\frak X(M)$ we get
$$\frac{\lb_0}{f^2}d\tau=-\frac12dZ+\frac12\chi dQ-\sigma_0d\tau,$$
and
$$dZ=\chi dQ-2\sigma d\tau\tag 3.5$$ where
$\sigma=\sigma_0+\frac{\lb_0}{f^2}$. From (3.5) we obtain
$$d\chi\w dQ-2d\sigma\w d\tau=0\tag 3.6$$
Since $d\chi=\gamma d\tau$ we have $d\tau\w(\gamma dQ+2d\sigma)=0$
which implies $\gamma dQ+2d\sigma=\kappa d\tau$ for a certain
function $\kappa$. Note that
$$\gather d\sigma_0=f^{-1}\a dQ+f^{-1}\phi dZ+(n-1)f^{-2}\phi^2
dQ+hd\tau=\\[-\frac1{4f^4}+(n-1)f^{-2}\phi^2]dQ+f^{-1}\phi dZ+hd\tau,
\endgather$$

for a certain function $h$. On the other hand
$\lb_0=-\frac{(n-4)Q}{4f^2}+\frac{C_0}{n-2}$. Hence
$d(\frac{\lb_0}{f^2})=-\frac{(n-4)}{4f^4}dQ+kd\tau$ and
$$d\sigma=[-\frac1{4f^4}+(n-1)f^{-2}\phi^2-\frac{(n-4)}{4f^4}]dQ+f^{-1}\phi
dZ+ld\tau$$  for some functions $k,l$. Since
$d\chi=(n-2)d(f^{-1}\phi)=-\frac{(n-2)}{4f^4}(1-4\phi^2f^2)d\tau$
we have $\gamma=-\frac{(n-2)}{4f^4}(1-4\phi^2f^2)$ and

$$[-\frac{3n-8}{4f^4}+\frac{(3n-4)\phi^2}{4f^2}]dQ+2f^{-1}\phi
dZ=md\tau,\tag 3.7$$ for a certain function $m$.

From equations (3.5) and (3.7) it follows that $$dQ\w d\tau=dZ\w
d\tau=0\tag 3.8$$  on a dense subset of $M$ and hence everywhere.

Define  a distribution $\Cal D=span\{\xi,J\xi\}$ and let $\Cal
D^{\perp}$ be an orthogonal (with respect to $g$ so also with
respect to $g_1$) complement to $\Cal D$. Both distributions are
defined in an open dense subset $U=\{x:\xi_x\ne0\}$.  Let
$\pi_{\De},\pi_{\DE}$ be orthogonal projections on $\De,\DE$
respectively.  Let us define $\om_{\De}(X,Y)=g_1(J\pi_{\De}X,Y),
\om_{\DE}(X,Y)=g_1(J\pi_{\DE}X,Y)$.  Then
$\om_{\De}+\om_{\DE}=\Om$ where $\Om(X,Y)=g_1(JX,Y)$ is the
K\"ahler form of $(M,g_1,J)$.  Note that
$\om_{\De}=\frac1{Q}d\tau\w d^c\tau$. Since $\xi$ is a holomorphic
Killing field on  $(M,g_1,J)$ it follows that
$H^{\tau}(JX,Y)=\frac12dd^c\tau(X,Y)$.  Since
$\n^1_vv=-\frac12\n^1Q=cv$ it follows that $\De$ is an
eigendistribution of both $\rho^1$ and $dd^c\tau$.  We have  (we
denote the Ricci form also by $\rho^1$)
$$\gather   \rho^1=\lb\om_{\De}+\om_1\tag
3.9a\\\frac12dd^c\tau=\mu\om_{\De}+\om_2,\tag 3.9b\endgather$$
where $\lb,\mu$ are eigenvalues of $\rho^1, H^{\tau}$
corresponding to an eigendistribution $\De$.  The eigenvalue $\mu$
satisfies an equation $\mu
Q=H^{\tau}(\n^1\tau,\n^1\tau)=-\frac12dQ(\n^1\tau)=-\frac12\beta
 Q$ where $dQ=\beta d\tau$.  Hence $\mu=-\frac12\beta$ and $d\mu\w
 d\tau=0$. From (3.4) it is clear that also $d\lb\w d\tau=0$.
Now we have
$$\widetilde{\rho}=\rho^1-\frac12\chi dd^c\tau-\sigma_0\Om,\tag 3.10$$
and consequently

$$\frac{\lb_0}{f^2}\om_{\De}=\lb\om_{\De}-\chi\mu\om_{\De}-\sigma_0\om_{\De}\tag
3.11$$ and
$$\frac{\lb_1}{f^2}\om_{\DE}=\om_1-\chi\om_2-\sigma_0\om_{\DE}.\tag
3.12$$

From (3.11) we obtain $\lb-\mu\chi=\sigma_0+\frac{\lb_0}{f^2}$.
Hence
$$\om_1-\chi\om_2=(\sigma_0+\frac{\lb_1}{f^2})\om_{\DE}=
(\lb-\mu\chi+\frac{\lb_1-\lb_0}{f^2})\om_{\DE}=\sigma_1\om_{\DE}.\tag 3.13$$

From (3.9) we get   $$\lb d\om_{\De}=-d\om_1, \mu
d\om_{\De}=-d\om_2.$$ Equation (3.13) implies that
$$
d\om_1-d\chi\w\om_2-\chi
d\om_2=d\sigma_1\w\om_{\DE}+\sigma_1d\om_{\DE},$$ thus
$$(-\lb+\mu\chi+\sigma_1)d\om_{\De}=d\chi\w\om_2+d\sigma_1\w\om_{\DE}.$$

Note that $d\om_{\De}=d(\frac1Qd\tau\w d^c\tau)=-\frac1Qd\tau\w
dd^c\tau=-\frac2Qd\tau\w (\mu\om_{\De}+\om_2)=-\frac2Qd\tau\w
\om_2$ and $\lb_1-\lb_0=\frac14(n-2)\frac Q{f^2}$. Let us write
$d\sigma_1=\psi d\tau$, then we obtain
$$d\tau\w(\frac{n-2}{4f^4}(-1-4f^2\phi^2)\om_2-\psi\om_{\DE})=0.\tag
3. 14$$  From (3.14) it is clear that in $U$ we have
$\om_2=\kappa_2\om_{\DE}$ for a certain function $\kappa_2\in C
^{\infty}(U)$. Hence also $\om_1=\kappa_1\om_{\DE}$ for a certain
function $\kappa_1\in C ^{\infty}(U)$.  It follows that the
function $\tau$ is a K\"ahler-Ricci potential in the sense  of
[D-M-1], [D-M-2], i.e. the distributions $\De,\DE$ are
eigendistribution of both $\rho^1$ and $H^{\tau}$. The fact that
the Einstein-K\"ahler manifold $(N,h)$ has a positive scalar
curvature is proved below. It is easy to check that also for $dim
M =4$ the manifold $(N,h)$ has constant scalar curvature. The
Einstein-Weyl structure on these manifolds is described in [W-W],
[M-P-P-S]. $\k$

\bigskip
{\bf 4. Eigenvalues of the Ricci tensor.} In our construction we
shall follow L. B\'erard Bergery (see [Ber],[J-2]).
 Let $(N,h,J)$ be a compact K\"ahler
Einstein manifold and $\dim N=2m$, $s\ge0,L>0,s\in \Bbb Q,L\in
\Bbb R$, and $g:[0,L]\rightarrow \Bbb R$ be a positive, smooth
function on $[0,L]$  which is even at $0$ and $L$, i.e. there
exists an $\e>0$ and  even, smooth functions
$g_1,g_2:(-\e,\e)\rightarrow \Bbb R$ such that $g(t)=g_1(t)$ for
$t\in[0,\e)$ and $g(t)=g_2(L-t)$ for $t\in(L-\e,L]$. Let
$f:(0,L)\rightarrow \Bbb{R}$ be positive on $(0,L)$, $f(0)=f(L)=0$
and let $f$ be odd at the points $0,L$. Let $P$ be a circle bundle
over $N$ classified by the integral cohomology class $\frac
s2c_1(N)\in H^2(N,\Bbb R)$ if $c_1(N)\ne 0$. Let $q$ be the unique
positive integer such that $c_1(N)=q\a$ where $\a\in H^2(N,\Bbb
R)$ is an indivisible integral class. Such a $q$ exists if $N$ is
simply connected or $\dim N=2$.  Note that every K\"ahler Einstein
manifold with positive scalar curvature is simply connected.  Then
$$s=\frac{2k}q; k\in\Bbb Z.$$ It is known that $q=n$ if
$N=\Bbb{CP}^{n-1}$ (see [B], p.273). Note that
$c_1(N)=\{\frac1{2\pi}\rho_N\}=\{\frac{\tau_N}{4m\pi}\0_N\}$ where
$\rho_N=\frac{\tau_N}{2m}\0_N$ is the Ricci form of $(N,h,J)$,
$\tau_N$ is the scalar curvature of $(N,h)$ and $\0_N$ is the
K\"ahler form of $(N,h,J)$. We can assume that $\tau_N=\pm 4m$. In
the case $c_1(N)=0$ we shall assume that $(N,h,J)$ is a Hodge
manifold, i.e. the cohomology class $\{\frac s{2\pi}\0_N\}$ is an
integral class. On the bundle $p:P\rightarrow N$ there exists a
connection form $\th$ such that $d\th=sp^*\0_N$ where
$p:P\rightarrow N$ is the bundle projection. Let us consider the
manifold $U_{s,f,g}=(0,L)\times P$ with the metric
$$k=dt^2+f(t)^2\th^2+g(t)^2p^*h.\tag 4.1$$
 It is known that the metric (4.1) extends to a metric
on the sphere bundle $M=P\times_{S^1}\Bbb{CP}^1$ if and only if a
function $g$ is positive and smooth on $[0,L]$, even at the points
$0,L$, the function $f$ is positive on $(0,L)$, smooth and odd at
$0,L$ and additionally
$$f'(0)=1,\qquad f'(L)=-1\tag 4.2$$ Then the metric (4.1) is
bi-Hermitian  (see [J-2]).  Note that $M=\Bbb P(L\oplus \Cal O) $
where $L=P\times_{S^1}\Bbb C$ with $S^1$ acting in a standard way
on $\Bbb C$ and $\Cal O$ is the trivial line bundle over $N$.

The metric $k=k_{f,g}$ extends to a metric on $\Bbb{CP}^n$ if and
only if the function $g$ is positive and smooth on $[0,L)$, even
at $0$, odd at  $L$, the function $f$ is positive, smooth and odd
at
 $0,L$ and additionally
$$f'(0)=1,\qquad f'(L)=-1,\qquad g(L)=0,\qquad g'(L)=-1.\tag 4.3$$
\bigskip
  Let us assume that
$(N,h)$ is a $2(n-1)$-dimensional  K\"ahler-Einstein manifold of
scalar curvature $4(n-1)\e$ where $\e\in\{-1,0,1\}$.  Using the
results in Sections 3 and [J-2]  we obtain the following formulae
for the eigenvalues of the Ricci tensor $\rho$ of
$(U_{s,f,g},k_{f,g})$:
$$\gather
\lb_0=-2(n-1)\frac{g''}g-\frac{f''}f,\tag 4.2\\
\lb_1=-\frac{f''}f+2(n-1)\bigg(\frac{s^2f^2}{4g^4}-\frac{f'g'}{fg}\bigg),\\
\lb_2=-\frac{g''}g+\bigg(\frac{s^2f^2}{4g^4}-\frac{f'g'}{fg}\bigg)+
\frac{2\e}{g^2}-\frac{3s^2f^2}{4g^4}-(2n-3)\frac{(g')^2}{g^2}.\endgather$$
We shall show that in fact $\e=1$, i.e. the scalar curvature of
the Einstein manifold $(N,h,J)$ is positive. From [J-1],
p.17,th.3.8. it follows that the conformal scalar curvature of
Einstein-Weyl manifold and hence $\lb_1$ is nonnegative.  We also
have for the Gauduchon metric $\lb_0=\lb_2$ and
$\lb_1+C^2f^2=\lb_0$ for a positive constant $C$.  Since
$f(0)=0=f(L)$ it follows that $f$ attains a maximum at a point
$t_0\in (0,L)$. Then $f'(t_0)=0$ and $f''(t_0)\le 0$. Hence at
$t_0$ we have
$$\lb_1=-\frac{f''}f+2(n-1)\frac{s^2f^2}{4g^4}>0$$
and
$$-\frac{g''}g+\frac{2\e}{g^2}-\frac{s^2f^2}{2g^4}-(2n-3)\frac{(g')^2}{g^2}=-2(n-1)\frac{g''}g-\frac{f''}f$$

and hence
$$\frac{2\e}{g^2}=\frac{s^2f^2}{2g^4}+(2n-3)\frac{(g')^2}{g^2}-(2n-3)\frac{g''}g-\frac{f''}f.$$

From (4.2) it follows that at $t_0$
$$-2(n-1)\frac{g''}g=2(n-1)\frac{s^2f^2}{4g^4}+C^2f^2>0$$ and
consequently  $\e>0$.
\medskip
{\it Acknowledgements.} The author thanks the referee for his
valuable remarks, which improved the paper.
\par
\bigskip

\centerline{\bf References.}
\par
\medskip
[B] A. Besse `Einstein manifolds' {\it Springer Verlag} Berlin
Heidelberg (1987)
\par
\medskip
\cite{Ber}  L. B\'erard Bergery,{\it Sur de nouvelles vari\'et\'es
riemanniennes d'Einstein}, Pu\-bl. de l'Institute E. Cartan
(Nancy) {\bf 4},(1982), 1-60.

\par
\medskip
[D-T]  M.  Dunajski, K.P. Tod  {Four dimensional metrics conformal
to K\"ahler}  to appear in Mathematical Proceedings of the
Cambridge Philosophical Society.

\par
\medskip
\cite{D-M-1} A. Derdzi\'nski, G. Maschler {\it Special
K\"ahler-Ricci potentials on compact K\" ahler manifolds},  J.
reine angew. Math. {\bf 593}, (2006), 73-116.
 \par
\medskip
\cite{D-M-2}  A. Derdzi\'nski, G.  Maschler {\it Local
classification of conformally-Einstein  K\"ahler metrics in higher
dimension},  Proc. London Math. Soc. (3) {\bf 87}, (2003), no. 3,
779-819.
\par
\medskip
\cite{D-M-3} A. Derdzinski , G. Maschler {\it A moduli curve for
compact conformally-Einstein K\"ahler manifolds},
 Compos. Math. 141 (2005), no. 4,1029-1080
 \par
\medskip
[G] P. Gauduchon {\it La 1-forme de torsion d'une variete
hermitienne compacte} Math. Ann. {\bf 267} (1984), 495-518.
\par
\medskip
[Gr] A. Gray {\it Einstein-like manifolds which are not Einstein}
Geom. Dedicata {\bf 7} (1978) 259-280.
\par
\medskip
[J-1] W. Jelonek {\it Killing tensors and Einstein-Weyl geometry}
Colloquium Math., vol.{\bf 81}, no.1, (1999), 5-19.
\par
\medskip
[J-2]W.Jelonek, {\it Higher dimensional Gray Hermitian manifolds},
J. Lond. Math. Soc. (2) 80 (2009), no. 3, 729-749.
\par
\medskip
[K] S. Kobayashi  {\it Transformation groups in differential
geometry}  Springer-Verlag 1972
\par
\medskip
[L] A. Lichn´erowicz, G´eom´etrie des groupes de transformations,
Dunod, Paris, 1958.
\par
\medskip
[M] A. Moroianu Lectures on K\"ahler Geometry, LMSST 69, Cambrige
University Press, 2007
\par
\medskip
[MPPS] B.Madsen, H. Pedersen, Y. Poon, A. Swann {\it Compact
Einstein-Weyl manifolds with large symmetry group.} Duke Math. J.
{\bf 88} (1997), 407-434.
\par
\medskip
[PS1] H. Pedersen and A. Swann {\it Riemannian submersions,
four-manifolds and Einstein-Weyl geometry} Proc.London Math. Soc
(3) {\bf 66} (1993) 381-399.
\par
\medskip
[PS2] H. Pedersen and A. Swann {\it Einstein-Weyl geometry, the
Bach tensor and conformal scalar curvature} J. reine angew.
Math.{\bf 441} (1993), 99-113.
\par
\medskip
[T]  K.P. Tod {\it Compact 3-dimensional Einstein-Weyl structures}
J.  London Math.  Soc.(2), 45 (1992) 341-351
\par
\medskip
[WW]  J. Wang, M.Y. Wang  {\it Einstein metrics on $S^2$-bundles}
Math. Ann. 310, 497-526, (1998)

\par
\medskip
\noindent Institute of Mathematics
\par
\noindent  Cracow University of Technology
\par
\noindent Warszawska 24
\par
\noindent 31-155 Krak\'ow, POLAND.
\par
\noindent e-mail: wjelon\@pk.edu.pl
\enddocument